\documentclass[leqno,letterpaper]{article}
\usepackage{latexsym,amsmath}
\usepackage{dsfont}
\long\def\ignore#1{}

\newtheorem{theorem}{Theorem}

\newtheorem{lemma}[theorem]{Lemma}

\newtheorem{coro}[theorem]{Corollary}

\def\qed{\hfill$\Box$\medskip}
\def\N{{\mathds N}}
\def\R{{\mathds R}}
\def\C{{\cal C}}
\def\F{{\cal F}}

\title{On distinct distances in homogeneous sets in the Euclidean space\thanks{The research was supported by OTKA and NSERC grants.}}

\author{
J\'ozsef Solymosi\thanks{Department of Mathematics, University of
British Columbia, Vancouver, BC, Canada V6T 1Z2,
\texttt{solymosi@math.ubc.ca}} \and Csaba D. T\'oth\thanks{Department of
Mathematics, MIT, Cambridge, MA~02139, USA,
\texttt{toth@math.mit.edu}}}

\date{}

\begin{document}
\maketitle

\begin{abstract}
 It is shown that every homogeneous set of $n$ points in the
 $d$-dimensional Euclidean space determines at least
 $\Omega(n^{2d/(d^2+1)} / \log^{c(d)} n)$ distinct distances for a
 constant $c(d)>0$. In three-space, the above general bound is
 slightly improved and it is shown that every homogeneous set of $n$
 points determines at least $\Omega(n^{.6091})$ distinct distances.
\end{abstract}

\section{Introduction\label{sec:intro}}

The history of the {\em distinct distance problem} goes back to
Erd\H os~\cite{e-osdnp-46} who asked the question: What is the
minimal number $g_d(n)$ of distinct distances determined by $n$
points in the $d$-dimensional Euclidean space $\R^d$? $n$ points in
the $d$-dimensional integer grid $[1,2,\ldots ,n^{\frac{1}{d}}]^d$
show that $g_d(n)= O(n^{2/d})$ for any $d\geq 2$ and, in particular,
$g_2(n)= O(n/\sqrt{\log n})$. Erd\H os conjectured that these bounds
are essentially optimal.

An initial lower bound of $g_2(n) \geq \Omega(\sqrt{n})$ by
Erd\H os~\cite{e-osdnp-46} was improved over the last almost $60$
years by Moser, Beck, Chung, Szemer\'edi, Trotter, and
Sz\'ekely~\cite{Moser,b-olpps-83,Chung,CST,s-cnhep-97}. Research
efforts have lead to several powerful methods (such as the crossing
theory~\cite{s-cnhep-97} and the $\varepsilon$-cutting
theory~\cite{ceg-ccbac-90}) which, in turn, found innumerable
applications in discrete and computational geometry. An excellent
survey by Pach and Sharir~\cite{ps-gi-04} elaborate on the history
of the distinct distance problem and its connections to other fields
of discrete mathematics. Determining the order of magnitude of
$g_2(n)$ (and $g_d(n)$ for every $d\in \N$) seems elusive. The
currently known best lower bound in the plane, $g_2(n) =
\Omega(n^{.8641})$, is due to Katz and Tardos~\cite{kt-neied-04}.
Their proof combines a method of Solymosi and
T\'oth~\cite{st-ddp-01} with results from entropy and additive
number theory.

In higher dimensions, not much work has been done. After some initial
results by Clarkson et al.~\cite{ceg-ccbac-90} and by Spencer et
al.~\cite{Spencer}, Aronov et al.~\cite{aps-ddthd-04} have showed
recently that the number of distinct distances determined by a set of
$n$ points in three-dimensional space is $g_3(n)=\Omega(n^{77/141
-\varepsilon}) =\Omega(n^{.5460})$ for any $\varepsilon> 0$. Solymosi
and Vu~\cite{sv-nobed-05} proved a general lower bound of
$g_d(n)=\Omega(n^{2/d-2/d(d+2)})$ for any fixed $d\geq 4$.

In this paper, we consider the minimum number $h_d(n)$ of distinct
distances in homogeneous sets of $n$ points in $\R^d$. A finite
point set $P\subset \R^d$ is {\em homogeneous} if the following two
conditions hold: $P$ lies in the interior of an axis-aligned
$d$-dimensional cube $C$ of volume $|P|$, and any unit cube in
$\R^d$ contains at most $O(1)$ points of $P$. Homogeneous sets
represent an important special case for the distinct distance
problem because the best known upper bound constructions (the
$d$-dimensional integer grids) are homogeneous, and because of
numerous connections to harmonic
analysis~\cite{b-sttsp-99,hi-cafed-05,ikp-fbedp-99,kt-scbfc-01}.
Iosevich~\cite{i-ccft-01} studied the distinct distance problem for
homogeneous sets (with additional restrictions). He showed that
$h_d(n)= \Omega(n^{3/2d})$ for any fixed $d\geq 2$.  Solymosi and
Vu~\cite{sv-ddhdh-04} proved a general bound of $h_d(n) =
\Omega(n^{2/d-1/d^2})$ for every dimension $d\geq 2$. For $d = 3$,
they have also obtained a slightly better bound $h_3(n) =
\Omega(n^{.5794})$.
In this paper, we improve all previous lower bounds on the number of
distinct distances in homogeneous sets of $n$ points in $\R^d$.

\begin{theorem}\label{thm:main}
For every $d\in \N$, there is a constant $c_d$ such that in every
homogeneous set $P$ of $n$ points in $\R^d$, there is a point $p\in
P$ from which there are at least
$$c_d \, n^{\frac{2d}{d^2+1}}
\log^{\frac{1-d^2}{d^2+1}} n$$
distinct distances measured to other points of $P$. In particular,
we have $h_d(n) \geq c_d \, n^{\frac{2d}{d^2+1}}
\log^{\frac{1-d^2}{d^2+1}} n$.
\end{theorem}

For $d=3,4$, and $5$, our general lower bound is $h_3(n) =
\Omega(n^{.5999})$, $h_4(n) = \Omega(n^{.4705})$, and $h_5(n) =
\Omega(n^{.3846})$. In three-dimensions, we slightly improve on this
bound and prove the following.

\begin{theorem}\label{thm:3D}
In every homogeneous set $P$ of $n$ points in $\R^3$, there is a
point $p\in P$ from which there are at least
$$\Omega \left(n^{\frac{53}{87}} \right) = \Omega(n^{.6091})$$
distinct distances measured to other points of $P$. In particular,
we have $h_3(n) = \Omega \left(n^{\frac{53}{87}} \right)$.
\end{theorem}

We prove Theorem~\ref{thm:main} in Section~\ref{sec:proof}. The proof
of Theorem~\ref{thm:3D} can be found in Section~\ref{sec:3D}.  In the
next section, we present a key lemma on the number of $k$-flats
incident to many points in a homogeneous point set in $\R^d$, for
$1\leq k < d$.

\section{Rich hyperplanes in homogeneous sets\label{sec:inc}}

Consider a set $P$ of $n$ points in $\R^d$. We say that a $k$-flat
(a $k$-dimensional affine subspace) is {\em $m$-rich} if it is
incident to at least $m$ points of $P$. The celebrated
Szemer\'edi-Trotter Theorem~\cite{st-epdg-83} states that for $n$
points in the plane, the number of $m$-rich lines ($1$-flats) is at
most $O(n^2/m^3+ n/m)$, and this bound is tight in the worst case.

The number of $m$-rich $k$-flats in $\R^d$ has been intensely
studied. The Szemer\'edi-Trotter type results have widespread
applications in discrete and combinatorial geometry. The
Szemer\'edi-Trotter Theorem's multi-dimensional
generalizations~\cite{aa-cfi-92,es-hipac-90,et-inh-05} always impose
some kind of restriction on the point set or on the set of $k$-flats,
otherwise $m$ points on a line give rise to an infinitely many of
$m$-rich $k$-flats for any $2\leq k\leq d$.

We adopt the following terminology: A set of $k+1$ points in $\R^d$,
$k\leq d$, is {\em affine independent} if it is contained in a
unique $k$-flat, which is said to be {\em span}ned by the point set.
A point set $P$ {\em determines} all the $k$-flats spanned by some
$k+1$ affine independent points of $P$.
For a constant $\alpha>0$, a finite point set $P\subset \R^d$ that
spans a $k$-flat is {\em $\alpha$-degenerate} if any $(k-1)$-flat
contains at most $\alpha \cdot |P|$ points of $P$.
For a finite point set $P\subset \R^d$ and a constant $\alpha>0$, we
say that a $k$-flat $F$ is {\em $\alpha$-degenerate} if the point set
$P\cap F$ is $\alpha$-degenerate. Note, for example, that all points
of $P\cap F$ in a 1-degenerate $k$-flat $F$ may lie on a $(k-1)$-flat,
but an $\alpha$-degenerate $k$-flat for $\alpha <1$ must be spanned by
points of $P$. We recall a result of Beck~\cite{b-olpps-83} on
$\alpha$-degenerate hyperplanes.

\begin{theorem}[Beck]\label{thm:Beck}
For every $k\in \N$, there are constants $\alpha_k,\beta_k>0$ with
the following property: For every $d\in \N$ and every finite point
set $P\subset \R^d$, if a $k$-flat $F$ is $\alpha_k$-degenerate,
then $P\cap F$ spans at least $\beta_k\cdot |F\cap P|^k$ distinct
$(k-1)$-flats.\qed
\end{theorem}

Elekes and T\'oth~\cite{et-inh-05} proved that for every dimension
$d\in \N$, there is a constant $\gamma_d>0$  such that the number of
$m$-rich $\gamma_d$-degenerate hyperplanes for $n$ points in $\R^d$
is at most $O(n^d/m^{d+1} + n^{d-1}/m^{d-1})$. The first term,
$O(n^d/m^{d+1})$, is dominant only if $m =O(\sqrt{n})$. We show
below a much stronger upper bound {\em for homogeneous sets}: A
homogeneous set of $n$ points in $\R^d$ determines at most
$O(n^d/m^{d+1})$ distinct $m$-rich hyperplanes for every $m\in \N$,
$d\leq m\leq n$.

We formulate our result for a slightly more general class of point
sets, where $n$ denotes the volume of the enclosing cube, rather than
the number of points. We say that a point set $P$ is {\em well
separated} if any unit cube in $\R^d$ contains at most $O(1)$ points
of $P$. By definition, every homogeneous set of $n$ points in $\R^d$
is well separated, and lies in a cube of volume $n$.

Let $f_{d,k}(P,m)$ denote the maximal number of $m$-rich $k$-flats in
a well separated point set $P$ contained in the interior of a
$d$-dimensional cube of volume $n$ in $\R^d$, and let
$$f_{d,k}(n,m) = \max_{P\subset \R^d, |P|=n} f_{d,k}(P,m).$$
Solymosi and Vu~\cite{sv-ddhdh-04} established the following lemma for
the number of $m$-rich lines in homogeneous sets of $n$ points in
$\R^d$. Their proof carries over verbatim for well separated sets of
volume $n$.

\begin{lemma}[Solymosi \& Vu]\label{lem:SV}
For every $d\in \N$, there is a constant $c_d$ such that
$$f_{d,1}(n,m)\leq c_d \, \frac{n^2}{m^{d+1}} .$$\qed
\end{lemma}

We extend their result for arbitrary $k\in \N$, $1\leq k\leq d-1$.

\begin{lemma}\label{lem:inc}
For every $d,k\in \N$, $1\leq k < d$, there is a constant $c_{d,k}$
such that
$$f_{d,k}(n,m)\leq c_{d,k} \, \frac{n^{k+1}}{m^{d+1}} .$$
\end{lemma}

The example of the $d$-dimensional integer grid $[1,2,\ldots
,n^{\frac{1}{d}}]^d$ shows that this bound is best possible for
every $m\in \N$, $1\leq m\leq n^{k/d}$.

\smallskip

\begin{proof}
For a fixed $d\in \N$, we prove that $f_{d,k}(n,m)=
O(n^{k+1}/m^{d+1})$. We proceed by induction on $k$, $1\leq k\leq
d$. The base case, $k=1$, is equivalent to Lemma~\ref{lem:SV}. Let us
assume that $1<k\leq d$ and that $f_{n_0,k_0}(P,m)=
O(n_0^{k_0+1}/m^{d+1})$ for every $k_0$, $1\leq k_0<k$, and $n_0\in
\N$.

Consider a well separated set $P$ that lies in the interior of a
$d$-dimensional cube $C$ of volume $n$. Clearly, we have
$|P|=O(n)$. We may choose an orthogonal coordinate system such that
all coordinates of every point of $P$ are irrational and $P$ lies in
the interior of cube $C$, whose vertices have rational
coordinates. This guarantees that for any subdivision of $C$ into
congruent subcubes, every point of $P$ lies in the {\em interior} of a
subcube.  For $i=0,1,\ldots ,\lceil\log n^{1/d}\rceil$, let $\C_i$
denote the subdivision of the cube $C$ into $2^{id}$ congruent
cubes. For instance, $\C_0=\{C\}$, $\C_1$ is a subdivision of $C$ into
$2^d$ cubes, and $\C_{\lceil (\log n)/d\rceil}$ is a subdivision into
constant volume cubes. There is a constant $\delta_d>d$ such that
every $k$-flat $F$ intersects at most $\delta_d |\C_i|^{k/d}=\delta_d
2^{ik}$ cubes of $\C_i$. If we put
$$\mu=\left\lfloor \frac{1}{k}\log
\frac{m}{4\delta_d(k+1)}\right\rfloor,$$
then every $m$-rich $k$-flat $F$ is incident to an average of at least
$m/(\delta_d 2^{\mu k}) \geq 4(k+1)$ points in a cube $Q\in
C_\mu$. That is, at least $m/2$ points of $P\cap F$ lie in subcubes
$Q\in C_\mu$ where $|P\cap F\cap Q|\geq 2(k+1)$.

Let $\alpha_k$ and $\beta_k$ be the constants from
Theorem~\ref{thm:Beck}. Let $\F$ denote the $m$-rich $k$-flats. We
classify the $k$-flats in $\F$ as follows:
\begin{itemize}
\item
    $\F_1=\{F\in \F: P\cap F$ is not $\alpha_k$-degenerate$\}$,
\item
    $\F_2=\{F\in \F:$ at least $\frac{m}{4}$ points of $P\cap F$
    lie in cubes $Q\in C_\mu$ such that the point set $P\cap
    F\cap Q$ is $\alpha_k$-degenerate$\}$,
\item  $\F_3=\F\setminus (\F_1\cup\F_2)$.
\end{itemize}
We show below that $|\F_q|=O(n^{k+1}/m^{d+1})$, for $q=1$, 2, and 3.
Every $F\in \F_1$ contains an $(\alpha_k m)$-rich $(k-1)$-flat. By
induction, the number of $(\alpha_k m)$-rich $(k-1)$-flats is
$O(n^k/(\alpha_k m)^{d+1}) = O(n^k/m^{d+1})$. Every $(\alpha_k
m)$-rich $(k-1)$-flat $R$ can be extended to a $m$-rich $k$-flat in
$O(n)$ different ways: $R$ together with a point of $P\setminus R$
spans a $k$-flat. This gives an upper bound
$|\F_1|=O(n^{k+1}/m^{d+1})$.

For an upper bound on $|\F_2|$, we consider the subdivision
$\C_\mu$. Let $K$ denote the affine independent $(k+1)$-element
subsets of $P$ that determine some $m$-rich $k$-flat in $\F_2$ and
lie in some cube $Q\in \C_\mu$. The volume of every cube $Q\in
\C_\mu$ is $O(n/2^{\mu d})=O(n/m^{d/k})$. Since $P$ is well separated,
we have $|P\cap Q|=O(n/m^{d/k})$. A trivial upper bound for the number
of affine independent $(k+1)$-element sets in all cubes of $\C_\mu$ is
$$|K|\leq |\C_\mu| \cdot
\left(O\left(\frac{n}{m^{d/k}}\right)\right)^{k+1}
=O\left(\frac{n^{k+1}}{m^d}\right).$$
We obtain a lower bound for $|K|$ by counting the affine independent
sets in each $F\in \F_2$. At least $m/4$ points of $P\cap F$ lie in
cubes $Q\in \C_\mu$ where the point set $P\cap F\cap Q$ is
$\alpha_k$-degenerate. By Theorem~\ref{thm:Beck}, every
$\alpha_k$-degenerate set $P\cap F\cap Q$ determines at least $\beta_k
|P\cap F\cap Q|^{k+1}$ affine independent $(k+1)$-element sets. If we
denote by $K(F)$ the number of $(k+1)$-element subsets of $K$ that
span $F$, then we have
$$|K(F)|\geq \mathop{\sum_{Q\in \C_{\mu}}}_{Q\cap F\neq
\emptyset}\beta_k |P\cap F\cap Q|^{k+1} \geq \delta_d
2^{\mu k} \left(\frac{m/4}{\delta_d 2^{\mu k}}\right)^{k+1} =
\Omega(m^{k+1} 2^{-\mu k^2}) = \Omega(m).$$

We conclude that $|K|=\sum_{F\in \F_2}=|\F_2|\cdot \Omega(m)$.  By
contrasting the upper and lower bounds for $|K|$, we get $|\F_2|=
O(n^{k+1}/m^{d+1})$.

Finally, we consider $\F_3$. For every $m$-rich $k$-flat $F\in \F_3$,
we define a set $S(F)$ of cubes from $\C_i$, $i=1,2,\ldots , \log
n^{1/d}$. A cube $Q\in \C_i$ is in $S(F)$ if and only if the point set
$P\cap F\cap Q$ is {\em not} $\alpha_k$-degenerate, but $P\cap F\cap
Q(i')$ is $\alpha_k$-degenerate for every $i'$, $0\leq i'<i$, where
$Q(i')$ is the (unique) cube $Q(i')\in C_{i'}$ containing $Q$.
If $P\cap F$ is not $\alpha_k$-degenerate, for example, then $C\not\in
S(F)$ for. Observe that the cubes of $S(F)$ are pairwise interior
disjoint and they jointly cover $P\cap F\cap C$.
We denote by $\mbox{dim}(X)$ the dimension of the affine subspace
spanned by a finite point set $X$.  For each $F\in F_3$, we further
classify the cubes in $S(F)$ according to three parameters:
For $i\in\{1,2,\ldots , \mu\}$, $j\in\{0, 1,\ldots , \log m\}$, and
$r\in \{1,\ldots , k-1\}$, let $S(F,i,j,r)$ denote the set of cubes
$Q\in S(F)$ such that
\begin{enumerate}
\item $Q\in \C_i$,
\item $2^{j-1}\cdot \frac{m}{\delta_d 2^{ik}}\leq |P\cap F\cap Q|
    < 2^j \cdot \frac{m}{\delta_d 2^{ik}}$,
\item $r=\min(k-1, \mbox{dim}(P\cap F\cap Q))$.
\end{enumerate}

Some of the cubes $Q\in S(F)$ are not included in any
$S(F,i,j,r,)\subset S(F)$: This is the case for every $Q\in S(F)\cap
\C_i$ for which $|P\cap F\cap Q|< (m/\delta_d 2^{ik+1})$ or $\mu<i$.
The total of number points of $P\cap F$ in these cubes is less than
$$
\mathop{\sum_{Q\in S(F)\cap \C_i}}_{0<i< \mu} |P\cap F\cap Q|
+ \mathop{\sum_{Q\in S(F)\cap \C_i}}_{i\geq \mu}
\frac{m}{\delta_d2^{ik}}< \frac{m}{2}+\frac{m}{4}= \frac{3m}{4}.$$
Therefore, the cubes in $S(F,i,j,r)$ for all $i$, $j$, $r$ jointly
contain at least $m/4$ points of $P\cap F$:
\begin{equation}\label{eq:m4}
\sum_{i=1}^{\mu} \sum_{j=0}^{\log m} \sum_{r=1}^{k-1}
|S(F,i,j,r)| \cdot m2^{j-ik}\geq \frac{m}{4}.
\end{equation}

For every $Q\in S(F,i,j,r)$, there is an $r$-flat $R\subset F$, such
that $|P\cap R\cap Q| \geq \alpha_k |P\cap F\cap Q| \geq \alpha_k
2^{j-1} m/(\delta_d 2^{ik}) = \Theta(2^{j-ik}m)$. Let us denote by
$Q'$ the cube in $\C_{i-1}$ that contains $Q\in \C_i$. Since $P\cap
F\cap Q'$ is already $\alpha_k$-degenerate, we have $|P\cap R\cap
Q|\leq \alpha_k |P\cap F\cap Q'|$. Let $D(Q,R)$ be the set of all
$(k-r)$-element affine independent sets $u\subset (P\cap F\cap
Q')\setminus R$ such that $R$ and $u$ together span $F$.  Since $P\cap
F\cap Q'$ is $\alpha_k$-degenerate, there are $\Theta(|P\cap F\cap
Q'|^{k-r})$ sets in $D(Q,R)$. Let $D'(Q,R)$ be a a subset of $D(Q,R)$
of size $\Theta(|P\cap F\cap Q|^{k-r})=\Theta((m2^{j-ik})^{k-r})$.

Let $T(F,i,j,r)$ denote the set of triples $(Q,R,u)$ such that $Q\in
S(F,i,j,r)$, $R$ is an $r$-flat with $|P\cap R\cap Q| \geq \alpha_k
|P\cap F\cap Q|$, and $u\in D'(Q,R)$. We have a lower bound
$$|T(F,i,j,r)|\geq |S(F,i,j,r)| \cdot \Theta((m2^{j-ik})^{k-r}).$$
Let us put
$$\tau(F,i,j,r)= \frac{|T(F,i,j,r)|}{(m2^{j-ik})^{k-r-1}},$$
and then Inequality~(\ref{eq:m4}) can be rewritten as
$$\sum_{i=1}^{\mu} \sum_{j=0}^{\log m} \sum_{r=1}^{k-1}
\tau(F,i,j,r) \geq \sum_{i=1}^{\mu} \sum_{j=0}^{\log
m} \sum_{r=1}^{k-1} |S(F,i,j,r)| \cdot \Omega(m2^{j-ik}) \geq
\Omega(m).$$
By summing over all $F\in \F_3$, we get
\begin{equation}\label{eq:fijr}
\sum_{F\in \F_3}\sum_{i=1}^{\mu} \sum_{j=0}^{\log m} \sum_{r=1}^{k-1}
\tau(F,i,j,r) \geq |\F_3|\cdot \Omega(m).
\end{equation}

We also compute an upper bound for the quantity on the left side of
Inequality (\ref{eq:fijr}). First, we give an upper bound on the
number of triples $(Q,R,u)\in T(F,i,j,r)$ for all $F\in \F_3$.  Recall
that $(Q,R,u)\in T(F,i,j,r)$ implies that $Q\in \C_i$, and $R$ is an
$r$-flat incident to $\ell=\Omega(m 2^{j-ik})$ points of $P\cap
Q$. Every cube $Q\in \C_i$ has volume $n/2^{ik}$ and $P\cap Q$ is well
separated. By our induction hypothesis, the number of $\ell$-rich
$r$-flats in $P\cap Q$ is $O((n/2^{ik})^{r+1}/\ell^{d+1})$. The cube
$Q'\in C_{i-1}$ contains $|P\cap Q'|=O(n/2^{(i-1)k})=O(n/2^{ik})$
points. So $P\cap Q'$ contains $(O(n/2^{ik}))^{k-r}$ distinct
$(k-r)$-element subsets. For all $Q\in \C_i$, we obtain an upper bound
$$\sum_{F\in \F_3} |T(F,i,j,r)| \leq |\C_i| \cdot O\left(
\frac{(n/2^{ik})^{r+1}}{(m2^{j-ik})^{d+1}}\right)
\cdot O\left(\left(\frac{n}{2^{ik}}\right)^{k-r}\right),$$
\begin{equation}\label{eq:T}
\sum_{F\in \F_3} |T(F,i,j,r)| \leq O\left( \frac{n^{k+1}}{m^{d+1}} \cdot 2^{ik-j(d+1)}\right).
\end{equation}
After dividing by $(m2^{j-ik})^{k-r-1}$, we sum
Inequality~(\ref{eq:T}) over all $i$, $j$, and $r$:
\begin{eqnarray}
\sum_{F\in \F_3} \tau(F,i,j,r) &\leq&
O\left( \frac{n^{k+1}}{m^{d+1}}  \cdot 2^{ik-j(d+1)}
\cdot \left(\frac{2^{ik}}{2^jm}\right)^{k-r-1}\right),\nonumber\\
\sum_{r=1}^{k-1}\sum_{F\in \F_3} \tau(F,i,j,r) &\leq&
O\left( \frac{n^{k+1}}{m^{d+1}}\cdot 2^{ik-j(d+1)}\right),\nonumber\\
\sum_{j=0}^{\log m}\sum_{r=1}^{k-1}\sum_{F\in \F_3} \tau(F,i,j,r)
&\leq& O\left( \frac{n^{k+1}}{m^{d+1}}\cdot 2^{ik}\right),\nonumber
\end{eqnarray}
\begin{equation}\label{eq:tau}
\sum_{i=1}^\mu\sum_{j=0}^{\log m}\sum_{r=1}^{k-1}\sum_{F\in \F_3}
\tau(F,i,j,r) \leq O\left( \frac{n^{k+1}}{m^{d+1}}\cdot m\right).
\end{equation}
By contrasting Inequalities~(\ref{eq:fijr}) and (\ref{eq:tau}), we
conclude that $|\F_3|= O(n^{k+1}/m^{d+1})$.  This completes the proof
of Lemma~\ref{lem:inc}.
\qed\end{proof}

\begin{coro}\label{cor:inc}
For every $d,k\in \N$, $1\leq k < d$, the number of incidences of
points and $m$-rich $k$-flats in a homogeneous set of $n$ points in
$\R^d$ is at most
$$O\left( \frac{n^{k+1}}{m^d} \right).$$
\end{coro}

\begin{proof}
In any homogeneous point set of size $n$ in $\R^d$, the number of
incidences of points and $m$-rich $k$-flats is bounded by
$$mf_{d,k}(P,m)+\sum_{j=m+1}^n f_{d,k}(P,j) \leq  O\left(
\frac{n^{k+1}}{m^d} \right)+ \sum_{j=m+1}^n O\left(
\frac{n^{k+1}}{j^{d+1}} \right) \leq O\left( \frac{n^{k+1}}{m^d}
\right).$$
\qed\end{proof}

\section{Proof of Theorem~\ref{thm:main}\label{sec:proof}}

We are given a homogeneous set $P$ of $n$ points in $d$-dimensions.
We may choose an orthogonal coordinate system such that all
coordinates of every point of $P$ are irrational and $P$ lies in the
interior of cube $C$, whose vertices have rational coordinates. This
guarantees that for any subdivision of $C$ into congruent subcubes,
every point of $P$ lies in the {\em interior} of a subcube.  Let $t$
denote the maximum number of distinct distances measured from a point
of $P$ (including distance $0$). There is a constant $\delta_d>d$ such
that for any $s\in \N$, any hyperplane or sphere intersects the
interior of at most $\delta_d s^{d-1}$ cubes in the subdivision of $C$
into $s^d$ congruent cubes.  We subdivide $C$ into $s^d$ congruent
subcubes $C_1,C_2,\ldots ,C_{s^d}$, where
$$s= \left\lfloor \left(\frac{n}{2\delta_d t} \right)^{\frac{1}{d-1}}
\right\rfloor.$$
Let $T$ be a set of triples $(p,q,c)\in P^3$ such that
\begin{itemize}
\item[(i)]   $p\neq q$,
\item[(ii)]  $p$ and $q$ lie in the same subcube $C_i$ for some $1\leq i \leq s^d$,
\item[(iii)] $p$ and $q$ are equidistant from $c$.
\end{itemize}

All points are located on $nt$ spheres centered at the $n$ points of
$P$.  The cubes $C_i$, $1\leq i\leq s^d$, subdivide each sphere into
{\em patches}. Since every sphere intersects at most $\delta_d
s^{d-1}$ subcubes $C_i$, there are at most $\delta_d n t
s^{d-1}=n^2/4$ patches, where each patch lies entirely in a subcube
$C_i$.  There are $n^2$ sphere-point incidences. The average number of
points on a patch is at least $4$. If $x$ points lie on a sphere patch
centered at $c$, then this patch contributes ${x\choose 2}2!$ triples
$(p,q,c)$ to $T$. We conclude that the number of triples is $|T| \geq
\Omega(n^2)$.

For every $m\in \N$, let $T_m$ denote the set of triples $(p,q,c)\in
T$ such that the bisector hyperplane of the segment $pq$ is incident
to at least $m$ but less then $2m$ points of $P$. Since every bisector
plane is incident to less than $n$ points, we can partition $T$ into
$\log n$ subsets
$$T = \bigcup_{j=0}^{\log n} T_{2^j}.$$
There is a value $m=2^j$ for some $0\leq j \leq \log n$, such that
$|T_m| \geq |T|/\log n\geq \Omega(n^2/\log n)$.

For a pair $(p,q)\in P^2$, $p\neq q$, all points of the set $M(p,q)=
\{c\in P : {\rm dist}(p,c) ={\rm dist}(q,c)\}$ lie on the bisector
hyperplane of the line segment $pq$. Every bisector hyperplane
intersects at most $\delta_d s^{d-1}$ subcubes, and in each subcube
$C_i$ it can bisect at most $|C_i\cap P|/2$ point pairs. So the number
of pairs $(p,q)\in P^2$ bisected by the same hyperplane is at most
$$\delta_d s^{d-1} \cdot O \left( \frac{n}{s^d}\right) = O \left(
\frac{n}{s}\right).$$

Let $B_m$ denote the set of all bisector hyperplanes that bisect the
pair $(p,q)$ for some $(p,q,c)\in T_m$. By definition, any hyperplane
in $B_m$ is incident to at least $m$ but less than $2m$ points of $P$.
By Lemma~\ref{lem:inc}, we have
$$|B_m| \leq O \left( \frac{n^d}{m^{d+1}} \right).$$

We can now give an upper bound for $|T_m|$. In a triple $(p,q,c)\in
T_m$, point $c$ lies on a bisector hyperplane of $B_m$. Each bisector
hyperplane is incident to less than $2m$ points of $P$ and bisects at
most $O(n/s)$ pairs $(p,q)$. Therefore
$$\Omega \left( \frac{n^2}{\log n} \right) \leq |T_m| \leq O \left(
\frac{n^d}{m^{d+1}} \right) \cdot 2m \cdot O \left( \frac{n}{s}
\right),$$
$$m^d \leq O\left( \frac{n^{d-1}\log n}{s} \right),
$$
\begin{equation}\label{eq:plug}
m\leq O \left( \frac{n^{\frac{d-1}{d}}\log^{1/d} n}{s^{1/d}} \right).
\end{equation}

We obtain another upper bound for $|T_m|$ by the following argument:
In a triple $(p,q,c)\in T_m$, both $p$ and $q$ lie in the same subcube
$C_i\subset C$. There are $s^d$ subcubes, and each subcube contains
$(O(n/s^d))^2 \leq O(n^2/s^{2d})$ point pairs. Hence, there are at
most $s^d \cdot O(n^2/s^{2d}) = O(n^2 / s^d)$ such pairs $(p,q)$. For
each pair $(p,q)$, where $(p,q,c)\in T_m$, there are at most $2m$
points $c\in P$ on the bisector hyperplane of $pq$. We conclude that
$$\Omega \left( \frac{n^2}{\log n} \right) \leq |T_m| \leq O\left(
\frac{n^2}{s^d} \right) \cdot 2m.$$
Using the upper bound for $m$ from Inequality (\ref{eq:plug}), we have
$$s^{\frac{d^2+1}{d}} \leq O \left( n^{\frac{d-1}{d}} \cdot
\log^{\frac{d+1}{d}} n\right),$$
$$\left( \frac{n}{t}\right)^{\frac{d^2+1}{d(d-1)}} \leq O \left(
n^{\frac{d-1}{d}} \cdot \log^{\frac{d+1}{d}} n\right),$$
$$\Omega \left( n^{\frac{2}{d-1}} \log^{-\frac{1+d}{d}} n\right) \leq
t^{\frac{d^2+1}{d(d-1)}},$$
$$\Omega \left(n^{\frac{2d}{d^2+1}} \log^{\frac{1-d^2}{d^2+1}} n
\right) \leq t,$$
as required. This completes the proof of Theorem~\ref{thm:main}
\qed

\section{Proof of Theorem~\ref{thm:3D}\label{sec:3D}}

Consider a homogeneous set $P$ of $n$ points in $\R^3$. Similarly to
the previous section, we assume that all coordinates of every point in
$P$ are irrational, and the vertices of the bounding cube $C$ have
rational coordinates. Let $t$ denote the maximum number of distinct
distances measured from a point of $P$ (including distance $0$). We
subdivide $C$ into $s^3$ congruent cubes $C_1,C_2,\ldots ,C_{s^3}$,
for
$$ s = \left\lfloor \sqrt{\frac{n}{\gamma t}} \right\rfloor,$$
where $\gamma>0$ is a constant to be specified later.

By Theorem~\ref{thm:Beck}, $P\subset \R^3$ contains $\Omega(n^2)$
affine independent point pairs. This implies that there is a subset
$P_0\subset P$ such that $|P_0|\geq \Omega(n)$ and every $c\in P_0$ is
incident to $\Omega(n)$ distinct lines spanned by $P$. For every $c\in
P_0$, let $P(c)\subset P\setminus \{c\}$ be a set of $\Omega(n)$
points such that the lines $cp$, $p\in P(c)$, are distinct.
For every point $c\in P_0$, let $H_c$ be a unit sphere centered at
$c$. For every $x\in \R^3\setminus \{c\}$, we denote by $\hat{x}$ the
projection of $x$ to the unit sphere $H_c$. Points of $P(c)$ have
distinct images in $H_c$ under this projection.  The set of images of
the projection is denoted by
$$\hat{P}(c):=\{ \hat{p} : c\in P(p)\}.$$

We partition the unit sphere $H_c$ into $6s^2$ convex spherical
regions $S_1(c),$ $S_2(c),\ldots ,S_{6s^2}(c)$ by $6s-12$ circular
arcs: Consider an axis-parallel cube centered at $c$ and subdivide
each of its 6 faces into $s^2$ congruent squares, then project these
squares to the sphere $H_c$ from $c$. The volume of each region is
$\Theta(1/s^2)$ and each region is contained in a disk of volume
$\Theta(1/s^2)$. Every circle on the sphere $H_c$ intersects at most
$O(s)$ regions. We then subdivide $\R^d\setminus \{c\}$ into $6s^2$
regions $R_i(c)$, $i=1,2,\ldots , 6s^2$, such that
$$R_i(c) = \{x \in \R^d\setminus \{c\}: \hat{x} \in S_i(c)\}.$$

For every $c\in P_0$ and $j=1,2,\ldots , 6s^2$, the region $R_j(c)$
contains $|P\cap R_i(c)|= O(n/s^2)$ points because the region
$R_j(c)\cap C$ can be covered by $O(n/s^2)$ unit cubes. Note also
that every plane incident to $c$ intersects at most $O(s)$ regions
$R_j(c)$, since every great circle of $S$ intersects at most $O(s)$
regions $S_j$. If $F$ is a plane, then $|F\cap R_j(c)\cap P| =
O(n^{2/3}/s)$ because $F\cap C$ can be covered by $O(n^{2/3})$ unit
cubes, and ${\rm area}(F\cap R_j(c)) \leq O({\rm area}(F\cap C)
/s)$.

For every $c\in P_0$, consider the at most $t$ spheres centered at $c$
that contain all points of $P(c)$. Every sphere $S$ centered at $c$ is
partitioned into {\em patch}es by the cubes $C_i$, $1\leq i\leq s^3$,
and the regions $R_j(c)$, $1\leq j\leq 6s^2$. We can partition $C$
into the subcubes $C_i$, $1\leq i\leq s^3$, by $3(s-1)$ planes.  These
planes partition every sphere $S$ along $3(s-1)$ circles.  Hence every
sphere $S$ is partitioned by $O(s)$ circular arcs into $O(s^2)$
patches. We partition the points of $P$ lying on a patch into disjoint
triples, after deleting at most two points from each patch if
necessary. This produces a set $Q$ of quadruples $(p,q,r,c)\in
P^3\times P_0$ such that,
\begin{itemize}
\item[(i)]   the points $p$, $q$, and $r$ are in $P(c)$;
\item[(ii)]  $p,q$, and $r$ lie on a sphere centered at $c$;

\item[(iv)] $p,q$, and $r$ lie in the same subcube $C_i$ for some $1\leq i \leq s^3$;
\item[(iv)] $p,q$, and $r$ lie in the same regions $R_j(c)$,
    for some $1\leq j\leq 6s^2$;
\item[(v)]   if $(p_1,q_1,r_1,c)\in Q$ and $(p_2,q_2,r_2,c)\in Q$, then
$\{p_1,q_1,r_1\} \cap \{p_2,q_2,r_2\}$  $= \emptyset$.
\end{itemize}

We give a lower bound on the number of quadruples in $Q$.
%
Let $g(c)$ denote the number of patches on all $O(t)$ spheres
centered at $c$: We have $g(c)=O(ts^2) = O(n/\gamma)$. The average
number of points on a patch centered at $c$ is $\Omega(\gamma n /
g(c)) = \Omega(\gamma)$. We choose the constant $\gamma>0$ such that
a patch contains at least $6$ points of $P(c)$ on average. If the
$k$-th patch contains a set of points $G_k(c)\subset P(c)$, then $Q$
contains $\lfloor |G_k(c)|/3\rfloor$ quadruples $(p,q,r,c)$. We
conclude that the total number of quadruples is
$$|Q| = \sum_{c\in P_0} \sum_{k=1}^{g(c)} \left\lfloor
\frac{|G_k|}{3} \right\rfloor \geq \Omega \left(n \sum_{k=1}^{g(c)}
\left( |G_k| - 2\right) \right)\geq\Omega(n^2).$$

We define the {\em multiplicity} of a pair $(p,q)\in P^2$ as
$$m(p,q)= |\{ c\in P_0: \exists r \mbox{ such that }
(p,q,r,c)\in Q\mbox{ or }(q,r,p,c)\in Q\mbox{ or }(r,p,q,c)\in
Q\}|.$$
We choose a parameter $m$ to be specified later, and distinguish two
types of quadruples in $Q$: A quadruple $(p,q,r,c)$ is {\em low} if at
least one edge of the triangle $pqr$ have multiplicity at most $m$.  A
quadruple $(p,q,r,c)$ is {\em high} if the multiplicity of all three
edges of $pqr$ are above $m$. Let $Q^-$ and $Q^+$ denote the sets of
low and high quadruples, respectively.
We distinguish two cases: First we consider the case that $|Q^+|\leq
|Q^-|$, then we proceed with the case $|Q^+|> |Q^-|$.

\smallskip

{\bf Case $|Q^+|\leq |Q^-|$.} There are at least $\Omega(n^2)$ low
quadruples in $Q$. We define a set of triples
$$ T :=  \{(p,q,c): (p,q,r,c)\in Q^-, m(p,q)\leq m\}.$$
We have extracted $|T|=\Omega(n^2)$ triples from $Q^-$. Similarly to
the previous section, we compute an upper bound on $|T|$. Every pair
$(p,q)$ from a triple of $T$ lies in one of the $s^3$ subcubes of
$C$, and for every pair $(p,q)$ there are at most $m$ centers $c$.
Therefore, we have an upper bound
$$|T| \leq s^3 \left( O\left( \frac{n}{s^3}\right) \right)^2 m =
O\left( \frac{mn^2}{s^3}\right).$$
Comparing this upper bound with the lower bound $|T|=\Omega(n^2)$,
we obtain
$$ \Omega(s^3)\leq m,$$
\begin{equation}\label{eq:m1}
\Omega \left( \frac{n^{3/2}}{t^{3/2}} \right) \leq m,
\end{equation}
$$\Omega \left( \frac{n}{m^{2/3}} \right) \leq t.$$

\smallskip

{\bf Case $|Q^+|> |Q^-|$.} At least half of the quadruples in $Q$
are high, and so $|Q^+|\geq \Omega(n^2)$.

For every $c\in P_0$, project the points of $P(c)$ to the sphere
$H_c$. If $(p,q,r,c)\in Q$, then the intersection of the bisector
plane of $pq$ and $H_c$ is the {\em bisector} (great circle) of the
segment $\hat{p}\hat{q}$ in the sphere $H_c$. A (possibly degenerate)
triangle $\hat{p}\hat{q}\hat{r}$ defines three distinct bisectors. The
bisectors of a triangle $\hat{p}\hat{q}\hat{r}$ meet in two antipodal
points on the sphere. The triangles that determine the same triple of
bisectors are similar (the center of similarity is the intersection of
the bisectors). Specifically, if the triangles
$\hat{p}_1\hat{q}_1\hat{r}_1,\hat{p}_2\hat{q}_2\hat{r}_2,
\ldots , \hat{p}_\ell\hat{q}_\ell\hat{r}_\ell$ determine the same
triple of bisectors, then the points $\hat{p}_1, \hat{p}_1, \ldots
,\hat{p}_\ell$ are collinear (the points $\hat{q}_1, \hat{q}_1,
\ldots ,\hat{q}_\ell$ and $\hat{r}_1, \hat{r}_1, \ldots
,\hat{r}_\ell$ are also collinear). Every triple of bisectors
determines a {\em family} of triangles. We define a {\em family of
quadruples} to be a collection of quadruples $(p,q,r,c)\in Q^+$ with a
common center $c$ such that the triangles $\hat{p}\hat{q}\hat{r}$ form
a family.

For every $c\in P_0$, we define a set of triangles in the sphere
$H_c$ by
$$ T(c) = \{\hat{p}\hat{q}\hat{r} : (p,q,r,c)\in Q^+\}.$$
By construction, all these triangles have pairwise disjoint vertex
sets. There is a set $P_1\subseteq P_0$ of size $\Omega(n)$ such that
for every $c\in P_1$, we have $|T(c)|=\Omega(n)$ triangles. For a
point $c\in P_1$, let $B_c$ denote the set of $m$-rich planes incident
to $c$. We denote by $\hat{B}_c$ the set of intersections of planes in
$B_c$ and the unit sphere $H_c$, which are great circles on $H_c$.
Note that the bisector of every edge $\hat{p}\hat{q}$ of a triangle of
$T(c)$ is in $\hat{B_c}$.

For $c\in P_1$, let us consider the partition of the sphere $H_c$ into
$6s^2$ regions $S_j(c)$, $1\leq j\leq 6s^2$, defined above.  Each
triangle of $T(c)$ lies entirely in one of the regions.  Let us denote
by $T_j(c)$ the set of triangles of $T(c)$ in $S_j(c)$ for every
$j=1,2,\ldots ,6s^2$. Since the triangles have disjoint vertex sets,
we have $|T_j(c)|\leq |P\cap R_j(c)|/3 \leq O(n/s^2) = O(t)$. But
$\sum_{j=1}^{6s^2} |T_j(c)|=\Omega(n)$, and so there are $\Omega(s^2)$
indices $j$ such that $|T_j(c)|= \Omega(n/s^2) = \Omega(t)$. Vertices
of similar triangles lie on three main circles. We have shown that
every region $R_j(c)$ contains at most $O(n^{2/3}/s) =
O(n^{1/6}t^{1/2})$ coplanar points.  Hence, there are at least
$\Omega(t^{1/2}/n^{1/6})$ families of triangles in $T_j(c)$. Since
each such family determines three distinct bisectors of $\hat{B}(c)$,
the triangles in $T_j(c)$ determine
$$\Omega \left( \left( \frac{t^{1/2}}{n^{1/6}} \right)^{1/3} \right) =
\Omega \left(\frac{t^{1/6}}{n^{1/18}} \right)$$
distinct bisectors in $\hat{B_c}$.  A bisector crosses at most $O(s)$
regions, and so we obtain the same bisector of $\hat{B_c}$ from at
most $O(s)$ regions. We conclude that the number of bisectors
determined by the $\Omega(n)$ triangles of $T(c)$ is
$$|B_c| \geq \frac{\Omega(s^2)}{O(s)} \cdot \Omega \left(
\frac{t^{1/6}}{n^{1/18}} \right) \geq \Omega \left(
\sqrt{\frac{n}{t}} \cdot \frac{t^{1/6}}{n^{1/18}} \right) = \Omega
\left( \frac{n^{4/9}}{t^{1/3}} \right).$$

Each of the $\Omega(n)$ points of $P_1$ is incident to $\Omega (
n^{4/9}/t^{1/3})$ distinct $m$-rich planes. This gives $\Omega (
n^{13/9}/t^{1/3})$ incidences on $m$-rich planes of $P$. By
Corollary~\ref{cor:inc}, we have
$$\Omega \left( \frac{n^{13/9}}{t^{1/3}} \right) \leq O
\left( \frac{n^3}{m^3} \right),$$
\begin{equation}\label{eq:m2}
m \leq O ( n^{14/27} t^{1/9}),
\end{equation}
$$ \Omega \left(\frac{m^{9}}{n^{14/3}} \right) \leq t.$$

In both cases, we have derived lower bounds for $t$ in terms of $n$
and $m$. We choose $m\in \N$ such that we obtain the same result in
both cases.  By comparing Inequalities (\ref{eq:m1}) and
(\ref{eq:m2}), we have
$$ \Omega \left( \frac{n^{3/2}}{t^{3/2}} \right) \leq m \leq O (
n^{14/27} t^{1/9}),$$
\begin{equation}\label{eq:t}
 \Omega \left( n^{53/87} \right)\leq t.
\end{equation}
The choice $m= n^{17/29} $ establishes Inequality (\ref{eq:t}) in
both cases. This completes the proof. \qed

\end{document}